\documentclass{amsart}

\usepackage{amssymb,xspace,amsmath,amsthm}

\usepackage{algorithmic}
\usepackage{algorithm}
\usepackage{graphics}
\usepackage{epsfig}
\usepackage{amsfonts}
\usepackage{color}
\usepackage{enumerate}

\newcommand{\LCP}{LCP\xspace}
\newcommand{\LCCP}{LC\ensuremath{\mathcal C}P\xspace}
\newcommand{\C}{\ensuremath{\mathcal C}\xspace}

\newcommand{\sep}{\ensuremath{\mathcal Sep}\xspace}
\newcommand{\CSP}{C\ensuremath{\mathcal Sep}P\xspace}
\newcommand{\LCSP}{LC\ensuremath{\mathcal Sep}P\xspace}
\newcommand{\NPC}{$NP$-complete\xspace}
\newcommand{\PB}[1]{\textsf{\scshape{#1}}}
\newcommand{\Gsol}{\ensuremath{G_\text{sol}}}

\newtheorem{definition}{Definition}
\newtheorem{example}{Example}
\newtheorem{theorem}{Theorem}
\newtheorem{proposition}{Proposition}
\newtheorem{problem}{Problem}
\newtheorem{corollary}{Corollary}
\newtheorem{lemma}{Lemma}
\newtheorem{property}{Property}
\newtheorem{claim}{Claim}

\begin{document}

\title{Longest Common Separable Pattern between Permutations}
%\titlerunning{Separable Pattern}
\author{Mathilde Bouvel}
\author{Dominique Rossin}
\address{
CNRS, Universit\'e Paris Diderot,
  Laboratoire d'Informatique Algorithmique: Fondements et Applications, 2 Place Jussieu, Case 7014,\\
F-75251 Paris Cedex 05, France}
\email{(mbouvel,rossin)@liafa.jussieu.fr}
\author{St\'ephane Vialette}
%\authorrunning{Mathilde Bouvel et al.}
%\tocauthor{Mathilde Bouvel (Universit\'e Paris Diderot and CNRS),
%  Dominique Rossin (Universit\'e Paris Diderot and CNRS),
%  St\'ephane Vialette (Universit\'e de Paris-Sud)}
\address{
Laboratoire de Recherche en Informatique (LRI), b\^at.490, Univ. Paris-Sud XI, \\ F-91405 Orsay cedex, France}
\email{Stephane.Vialette@lri.fr}

\maketitle
\begin{abstract}
In this article, we study the problem of finding the longest common separable pattern between several permutations. We give a polynomial-time algorithm when the number of input permutations is fixed and show that the problem is NP-hard for an arbitrary number of input permutations even if these permutations are separable.

On the other hand, we show that the NP-hard problem of finding the longest common pattern between two permutations cannot be approximated better than within a ratio of $\sqrt{Opt}$ (where $Opt$ is the size of an optimal solution) when taking common patterns belonging to pattern-avoiding classes of permutations.
\end{abstract}

\section{Introduction and basic definitions}
\label{section:intro}

A permutation $\pi$ is said to be a pattern within a permutation
$\sigma$ if $\sigma$ has a subsequence that is order-isomorphic to
$\pi$. 
Within the last few years, the study of the \emph{pattern containment}
relation on permutations has become a very active area of research in
both combinatorics and computer science. 
In combinatorics, much research focused on closed classes of
permutations, \emph{i.e.}, permutations that are closed downwards under
forming subpermutations. 
A huge literature is devoted to this subject.
To cite only a few of a plethora of suitable examples,
Knuth considered permutations that do not contain the pattern
$3 1 2$~\cite{Knuth:ArtComputerProgramming:1:1973},
Lov\`asz considered permutations that do not contain the pattern
$2 1 3$~\cite{Lovasz:1979} and
Rotem those that do not contain $2 3 1$ nor 
$3 1 2$~\cite{Rotem:DM:1981}.  

Surprisingly enough, there is considerably less research on
algorithmic aspects of pattern involvement.
Actually, it appears to be a difficult problem to decide whether a
permutation occurs as a pattern in another permutation. Indeed, the problem in this general version  is $NP$-complete \cite{BBL98}. 
The case of \emph{separable patterns},
\emph{i.e.}, permutations that contain neither the subpattern
$3 1 4 2$ nor $2 4 1 3$, was proved to be
solvable in $\mathcal{O}(kn^6)$ time and
$\mathcal{O}(kn^4)$ space in \cite{BBL98},
where $k$ is the length of the pattern and $n$ is the length of the
target permutation.
The design of efficient algorithms for the recognition of a fixed
pattern in a permutation is considered in
\cite{Albert:Aldred:Atkinson:Holton:ISAAC:2001} and in particular a 
$\mathcal{O}(n^5 \log n)$ time algorithm is given for finding
separable patterns.  
L.~Ibarra subsequently improved the complexity for separable patterns
to $\mathcal{O}(kn^4)$ time and $\mathcal{O}(kn^3)$ space in
\cite{Iba97}.
Beside separable patterns, only a few restricted cases were
considered. 
A $\mathcal{O}(n \log \log n)$ time algorithm is presented in
\cite{Chang:Wang:IPL:1992} for finding the longest increasing or
decreasing subpermutation of a permutation of length $n$.

In the present paper we continue this line of research on separable
patterns by considering the problem of finding a maximum length common
separable pattern to a set of permutations, \emph{i.e.},
given a set of permutations, find a longest separable permutation
that occurs as a pattern in each input permutation.
Of particular importance in this context, we do not impose here the
input permutations to be separable.

This paper is organized as follows.
In the remainder of Section~\ref{section:intro},
we briefly discuss basic notations and definitions that we will use
throughout.
In Section~\ref{section:poly}, we give a polynomial-time algorithm
for finding the largest common separable pattern that appears as a
pattern in a fixed number of permutations. 
Section~\ref{section:hardness result} is devoted to proving hardness
of the problem.
Finally, some inapproximation issues are presented in
Section~\ref{section:approximation}.

\subsection{Permutations}

A permutation $\sigma \in S_n$  is a bijective map from $[1..n]$ to itself. The integer $n$ is called the \emph{length} of $\sigma$. We denote by $\sigma_i$ the image of $i$ under $\sigma$. A permutation can be seen as a word $\sigma_1 \sigma_2 \ldots \sigma_n$ containing exactly once each letter $i \in [1..n]$. For each entry $\sigma_i$ of a permutation $\sigma$, we call $i$ its \emph{index} and $\sigma_i$ its \emph{value}.

\begin{definition}[Pattern in a permutation]
A permutation $\pi \in S_k$ is a \emph{pattern} of a permutation $\sigma \in S_n$ if there is a subsequence of $\sigma$ which is order-isomorphic to $\pi$; in other words, if there is a subsequence $\sigma_{i_1} \sigma_{i_2} \ldots \sigma_{i_k}$ of $\sigma$ (with $1 \leq i_1 < i_2 <\ldots<i_k \leq n$) such that $\sigma_{i_{\ell}} < \sigma_{i_m}$ whenever $\pi_{\ell} < \pi_{m}$. \\
We also say that $\pi$ is \emph{involved} in $\sigma$ and call $\sigma_{i_1} \sigma_{i_2} \ldots \sigma_{i_k}$ an \emph{occurrence} of $\pi$ in $\sigma$. \label{def:pattern}
% whose entries appear in the same relative order as the entries of $\pi$ i.e. if there exist integers $1 \leq i_1 < i_2 <\ldots<i_k \leq n$ such that $\sigma_{i_{\ell}} < \sigma_{i_m}$ whenever $\pi_{\ell} < \pi_{m}$. \label{def:pattern}
\end{definition}

A permutation $\sigma$ that does not contain $\pi$ as a pattern is said to \emph{avoid} $\pi$. Classes of permutations of interest are the \emph{pattern-avoiding classes of permutations}: the class of all permutations avoiding the patterns $\pi_1, \pi_2 \ldots \pi_k$ is denoted $S(\pi_1, \pi_2, \ldots, \pi_k)$, and $S_n(\pi_1, \pi_2, \ldots, \pi_k)$ denotes the set of permutations of length $n$ avoiding $\pi_1, \pi_2, \ldots, \pi_k$.

\begin{example}
For example $\sigma=1  4  2  5  6  3$ contains the pattern $1  3  4  2$, and $1  5  6  3$, $1  4  6  3$, $2 5 6 3$ and $1 4 5 3$ are the occurrences of this pattern in $\sigma$. But $\sigma \in S(3  2  1)$: $\sigma$ avoids the pattern $3  2  1$ as no subsequence of length $3$ of $\sigma$  is isomorphic to $3  2  1$, \textit{i.e.}, is decreasing. \label{ex:pattern}
\end{example}

If a pattern $\pi$ has an occurrence $\sigma_{i_1} \sigma_{i_2} \ldots \sigma_{i_k}$ in a permutation $\sigma$ of length $n$, let $I$ and $V$ be two subintervals of $[1..n]$ such that $\{i_1, i_2, \ldots, i_k\} \subseteq I$ and $\{\sigma_{i_1}, \sigma_{i_2}, \ldots, \sigma_{i_k} \} \subseteq V$;  then we say that $\pi$ has an occurrence in $\sigma$ in the intervals $I$ of indices and $V$ of values, or that $\pi$ is a pattern of $\sigma$ using the intervals $I$ of indices and $V$ of values in $\sigma$.

Among the pattern-avoiding classes of permutations, we are particularly interested here in the separable permutations.
\begin{definition}[Separable permutation]
The class of separable permutations, denoted \sep, is $\sep = S(2413,3142)$.
\end{definition}

There are numerous characterizations of separable permutations, for
example in terms of permutation graphs \cite{BBL98}, of interval
decomposition \cite{BXHC05,BCMR05,br06}, or with ad-hoc structures
like the separating trees \cite{BBL98,Iba97}. Separable permutations
have been widely studied in the last decade, both from a combinatorial
\cite{Wes95,EHPR98} and an algorithmic \cite{BBCP07,BBL98,Iba97} point
of view. 

We define two operations of concatenation on permutation patterns:

\begin{definition}[Pattern concatenation]
Consider two patterns $\pi$ and $\pi'$ of respective lengths $k$ and $k'$. The positive and negative concatenations of $\pi$ and $\pi'$ are defined respectively by:\\
\begin{align*}
\pi \oplus \pi' &= \pi_1 \cdots \pi_k (\pi'_1 +k) \cdots (\pi'_{k'} +k) \\
\pi \ominus \pi' &= (\pi_1+k') \cdots (\pi_k+k') \pi'_1 \cdots \pi'_{k'} 
\end{align*}
 \label{def:concatenation}
\end{definition}
The following property, whose proof is straightforward with separating trees, is worth noticing for our purpose:
\begin{property}
If both $\pi$ and $\pi'$ are separable patterns, then $\pi \oplus \pi'$ and $\pi \ominus \pi'$ are also separable. Conversely, any separable pattern $\pi$ of length at least $2$ can be decomposed into $\pi=\pi_1 \oplus \pi_2$ or $\pi=\pi_1 \ominus \pi_2$ for some smaller but non-empty separable patterns $\pi_1$ and $\pi_2$. \label{prop:concatenation}
\end{property}

\subsection{Pattern problems on permutations}

The first investigated problem on patterns in permutations is the \emph{Pattern Involvement} Problem:
\begin{problem}[Pattern Involvement Problem]~\\
\textsc{Input}: A pattern $\pi$ and a permutation $\sigma$. \\ 
\textsc{Output}: A boolean indicating whether $\pi$ is involved in $\sigma$ or not.
\label{problem:involvement}
\end{problem}
It was shown to be $NP$-complete in \cite{BBL98}. However, in \cite{BBL98} the authors also exhibit a particular case in which it is polynomial-time solvable: namely when the pattern $\pi$ in input is a separable pattern.

Another problem of interest is the \emph{Longest Common Pattern} Problem (\PB{LCP} for short): 
\begin{problem}[\PB{LCP} Problem]~\\
\textsc{Input}: A set $X$ of permutations.\\
\textsc{Output}: A pattern of maximal length occurring in each $\sigma \in X$.
\label{problem:LCP}
\end{problem}

This problem is clearly $NP$-hard in view of the complexity of Problem \ref{problem:involvement}. We showed in \cite{br06} that it is solvable in polynomial time when $X= \{\sigma_1,\sigma_2\}$ with $\sigma_1$ a separable permutation (or more generally, when the length of the longest \emph{simple permutation} \cite{BHV06a,BHV06b,BRV06} involved in $\sigma_1$ is bounded).

In this paper, we will consider a restriction of Problem \ref{problem:LCP}. For any (pattern-avoiding) class \C of permutations, we define the \emph{Longest Common \C-Pattern} Problem (\PB{\LCCP} for short): 
\begin{problem}[\PB{\LCCP} Problem]~\\
\textsc{Input}: A set $X$ of permutations.\\
\textsc{Output}: A pattern of \C of maximal length occurring in each $\sigma \in X$.
\label{problem:LCCP}
\label{problem:LCSP}
\end{problem}

In particular, we focus in this paper on the \emph{Longest Common Separable Pattern} Problem (\PB{\LCSP}) which in fact is  \PB{\LCCP} where $\C = \sep$.

%% \begin{problem}[\PB{\LCSP} Problem]\\
%% \textsc{Input}: A set $X$ of permutations.\\
%% \textsc{Output}: A separable pattern of maximal length occurring in each $\sigma \in X$.
%% \label{problem:LCSP}
%% \end{problem}

To our knowledge,  complexity issues of  the \PB{\LCCP} Problem are still unexplored. We will show in this paper that the \PB{\LCSP} Problem is $NP$-hard in general, but solvable in polynomial-time when the cardinality of the set $X$ of permutations in input is bounded by any constant $K$. However the method we use in our polynomial-time algorithm for solving \PB{\LCSP} on $K$ permutations is specific to separable patterns and cannot be extended to any class \C of pattern-avoiding permutations.

Some classes \C of permutations are known for which even the \emph{Recognition} Problem (\textit{i.e.}, deciding if a permutation belongs to \C) is $NP$-hard, so that \PB{\LCCP} on $K$ permutations must be $NP$-hard for those classes. \cite{AAADHHO03} gives the example of the class of $4$-stack sortable permutations.

However, we are not aware of any  example of \emph{finitely based} pattern-avoiding classes of permutations (with a finite number of excluded patterns) for which the Recognition Problem is $NP$-hard. Thus an open question is to know if the \PB{\LCCP} problem for $K$ permutations is polynomial-time solvable for any finitely based \C, or to exhibit such a class \C for which this problem is $NP$-hard.

\section{Polynomial algorithm for the longest common separable pattern
  between a finite number of permutations}
\label{section:poly}

In \cite{BBL98}, the authors show that the problem of deciding whether a permutation $\pi$ is a pattern of a permutation $\sigma$ is $NP$-complete. A consequence is that the problem of finding a longest common pattern between two or more permutations in $NP$-hard. However, they describe a polynomial-time algorithm for solving the Pattern Involvement Problem when the pattern $\pi$ is \emph{separable}. This algorithm uses dynamic programming, and processes the permutation according to one of its separating trees.
%% following the structure of a binary separating tree of the separable pattern $\pi$.

With the same ideas, we described in \cite{br06} a polynomial-time algorithm for finding a longest common pattern between two permutations, provided that one of them is separable. Notice that a longest common pattern between two permutations, one of them being separable, is always separable.

In this section, we generalize the result obtained in \cite{br06} giving a polynomial-time algorithm for finding a longest common \emph{separable} pattern between $K$ permutations, $K$ being any fixed integer, $K \geq 1$. Notice that we make no hypothesis on the $K$ input permutations.

Like in \cite{BBL98} and \cite{br06}, our algorithm will use dynamic programming. However, since we do not have a separability hypothesis on any of the permutations, we cannot design an algorithm based on a separating tree associated to one of the permutations in input. To compute a longest common separable pattern between the input permutations, we will only consider sub-problems corresponding to $K$-tuples of intervals of indices and values, one such pair of intervals for each permutation.

Namely, let us consider $K$ permutations $\sigma^1, \ldots, \sigma^K$, of length $n_1, \ldots , n_K$ respectively, and denote by $n$ the maximum of the $n_q$'s, $1\leq q\leq K$. For computing a longest common separable pattern between $\sigma^1, \ldots , \sigma^K$, we will consider a dynamic programming array $M$ of dimension $4K$, and when our procedure for filling in $M$ ends, we intend that $M(i_1, j_1, a_1, b_1, \ldots, i_K, j_K, a_K, b_K)$ contains a common separable pattern $\pi$ between $\sigma^1, \ldots, \sigma^K$ that is of maximal length among those using, for any $q \in [1..K]$, intervals $[i_q..j_q]$ of indices and $[a_q..b_q]$ of values in $\sigma^q$. If we are able to fill in $M$ in polynomial time, with the above property being satisfied, the entry $M(1, n_1, 1, n_1, \ldots, 1, n_K, 1, n_K)$ will contain, at the end of the procedure, a longest common separable pattern between $\sigma^1, \ldots, \sigma^K$.\\

\par
Algorithm \ref{alg:separable} shows how the array $M$ can indeed be filled in in polynomial time. In Algorithm \ref{alg:separable}, $Longest$ is the naive linear-time procedure that runs through a set $S$ of patterns and returns a pattern in $S$ of maximal length.\\

\begin{algorithm}[!ht]
\caption{Longest common separable pattern between $K$ permutations}
\label{alg:separable}
\begin{algorithmic}[1]
\STATE {\textsc{Input}: $K$ permutations $\sigma^1, \ldots, \sigma^K$ of length $n_1, \ldots , n_K$ respectively \\~\\}
\STATE {\textsc{Create an array $M$:} }
\FOR{any integers $i_q$, $j_q$, $a_q$ and $b_q \in [1..n_q] $, for all $q \in [1..K]$}
\STATE {$M(i_1,j_1,a_1,b_1, \ldots,i_K,j_K,a_K,b_K) \leftarrow \epsilon$ }
\ENDFOR \\~\\
\STATE {\textsc{Fill in $M$:} }
\FOR{any integers $i_q$, $j_q$, $a_q$ and $b_q \in [1..n_q]$, $i_q\leq j_q$, $a_q \leq b_q$, for all $q\in [1..K]$, by increasing values of $\sum_q (j_q-i_q) +(b_q-a_q)$}
\IF {$\exists q \in [1..K]$ such that $i_q=j_q$ or $a_q=b_q$}
\IF {$\forall q \in [1..K], \exists h_q \in [i_q..j_q]$ such that $\sigma^q_{h_q} \in [a_q..b_q]$}
\STATE {$M(i_1,j_1,a_1,b_1, \ldots ,i_K,j_K,a_K,b_K) \leftarrow 1 $}
\ELSE 
\STATE {$M(i_1,j_1,a_1,b_1, \ldots ,i_K,j_K,a_K,b_K) \leftarrow \epsilon $}
\ENDIF
\ELSE
\STATE { /*$\forall q \in [1..K], i_q < j_q$ and $a_q<b_q$ /* \\ $M(i_1,j_1,a_1,b_1, \ldots ,i_K,j_K,a_K,b_K) \leftarrow Longest( S_{\oplus} \cup S_{\ominus} \cup S )$ where 
\begin{eqnarray*}
 S_{\oplus} & = & \{ M(i_1,h_1-1,a_1,c_1-1, \ldots ,i_K,h_K -1,a_K,c_K-1)  \oplus  M(h_1,j_1,c_1,b_1, \\ & & \ldots ,h_K,j_K,c_K,b_K) : i_q < h_q \leq j_q , a_q < c_q \leq b_q , \forall q \in [1..K] \} \\
 S_{\ominus} & = & \{ M(i_1,h_1-1,c_1,b_1, \ldots ,i_K,h_K -1,c_K,b_K)  \ominus  M(h_1,j_1,a_1,c_1-1, \\ & & \ldots ,h_K,j_K,a_K,c_K-1) : i_q < h_q \leq j_q , a_q < c_q \leq b_q , \forall q \in [1..K] \} \\
S_{~} & = & \{1\} \textrm{ if } \forall q \in [1..K], \exists h_q \in [i_q..j_q] \textrm{ such that } \sigma^q_{h_q} \in [a_q..b_q],\\
 & = &  \{\epsilon\} \textrm{ otherwise.}
\end{eqnarray*}
}
\ENDIF
\ENDFOR \\~\\
\STATE {\textsc{Output}: $M(1,n_1,1,n_1, \ldots, 1, n_K, 1, n_K)$}
\end{algorithmic}
\end{algorithm}

\par
Before giving the details of the proof of our algorithm for finding a longest common separable pattern, we state and prove two lemmas. They should help understanding how common separable patterns can be merged, or on the contrary split up, to exhibit other common separable patterns. We are also interested in the stability of the maximal length property when splitting up patterns.

\begin{lemma}
Let $\pi_1$ be a common separable pattern between $\sigma^1, \ldots, \sigma^K$ that uses the intervals $[i_q..h_q-1]$ of indices and $[a_q..c_q-1]$ (resp. $[c_q..b_q]$) of values in $\sigma^q$, for all $q \in [1..K]$. 

Let $\pi_2$ be a common separable pattern between $\sigma^1, \ldots, \sigma^K$ that uses the intervals $[h_q..j_q]$ of indices and $[c_q..b_q]$ (resp. $[a_q..c_q-1]$) of values in $\sigma^q$, for all $q \in [1..K]$. 

Then $ \pi = \pi_1 \oplus \pi_2$ (resp. $\pi = \pi_1 \ominus \pi_2$) is a common separable pattern between $\sigma^1, \ldots, \sigma^K$ that uses the intervals $[i_q..j_q]$ of indices and $[a_q..b_q]$ of values in $\sigma^q$, for all $q \in [1..K]$. \label{lem:fusion}
\end{lemma}

\begin{proof}
We give a proof for $\pi = \pi_1 \oplus \pi_2$ (the case $\pi=\pi_1\ominus\pi_2$ being similar). 

Fix some $q \in [1..K]$. By hypothesis, there exist occurrences of $\pi_1$ and $\pi_2$ in $\sigma^q$, the occurrence of $\pi_1$ using the intervals $[i_q..h_q-1]$ of indices and $[a_q..c_q-1]$ of values, and the occurrence of $\pi_2$ using the intervals $[h_q..j_q]$ of indices and $[c_q..b_q]$ of values. It is then easily noticed (see Figure \ref{fig:fusion}) that all the elements used in one of these occurrences form an occurrence of the pattern $\pi = \pi_1 \oplus \pi_2$ in $\sigma^q$ in the intervals $[i_q..j_q]$ of indices and $[a_q..b_q]$ of values. This argument holds for any $q \in [1..K]$ and hence $\pi$ is a common separable pattern between $\sigma^1, \ldots, \sigma^K$ using the intervals $[i_q..j_q]$ of indices and $[a_q..b_q]$ of values in $\sigma^q$, for all $q \in [1..K]$.

\begin{center}
\begin{figure}
\input{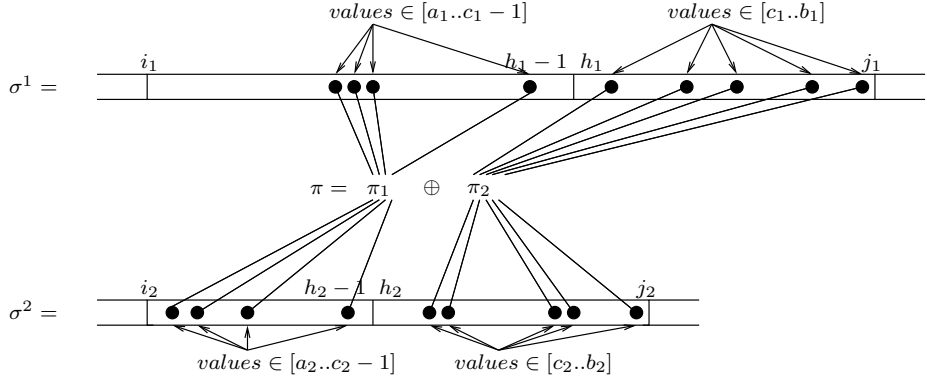}
\caption{Proof of lemma \ref{lem:fusion} for $K=2$}
\label{fig:fusion}
\end{figure}
\end{center}
\qed
\end{proof}

\begin{lemma} 
Let $\pi$ be a common separable pattern of maximal length between $\sigma^1, \ldots, \sigma^K$ among those using the intervals $[i_q..j_q]$ of indices and $[a_q..b_q]$ of values in $\sigma^q$, for all $q \in [1..K]$.

If $\pi = \pi_1 \oplus \pi_2$ (resp. $\pi= \pi_1 \ominus \pi_2$), with $\pi_1$ and $\pi_2$ \emph{non-empty} separable patterns, then there exist indices $(h_q)_{q \in [1..K]}$ and values $(c_q)_{q \in [1..K]}$, with $i_q < h_q \leq j_q , a_q < c_q \leq b_q , \forall q \in [1..K]$, such that: 
\begin{enumerate}[$i)$]
 \item $\pi_1$ is a common separable pattern of maximal length between $\sigma^1, \ldots, \sigma^K$ among those using the intervals $[i_q..h_q-1]$ of indices and $[a_q..c_q-1]$ (resp. $[c_q..b_q]$) of values in $\sigma^q$, for all $q \in [1..K]$, and
 \item $\pi_2$ is a common separable pattern of maximal length between $\sigma^1, \ldots, \sigma^K$ among those using the intervals $[h_q..j_q]$ of indices and $[c_q..b_q]$ (resp. $[a_q..c_q-1]$) of values in $\sigma^q$, for all $q \in [1..K]$.
\end{enumerate}
% $\pi_1$ is a common separable pattern of maximal length between $\sigma^1, \ldots, \sigma^K$ among those using the intervals $[i_q..h_q-1]$ of indices and $[a_q..c_q-1]$ (resp. $[c_q..b_q]$) of values in $\sigma^q$, for all $q \in [1..K]$ and $\pi_2$ is a common separable pattern of maximal length between $\sigma^1, \ldots, \sigma^K$ among those using the intervals $[h_q..j_q]$ of indices and $[c_q..b_q]$ (resp. $[a_q..c_q-1]$) of values in $\sigma^q$, for all $q \in [1..K]$.
\label{lem:cutting}
\end{lemma}

\begin{proof}
Again, consider the case $\pi = \pi_1 \oplus \pi_2$ (the case $\pi = \pi_1 \ominus \pi_2$ being similar). 

Fix some $q \in [1..K]$. By hypothesis, $\pi = \pi_1 \oplus \pi_2$ has an occurrence in $\sigma^q$ in the intervals $[i_q..j_q]$ of indices and $[a_q..b_q]$ of values. By definition of positive pattern concatenation, this occurrence splits into two occurrences of $\pi_1$ and $\pi_2$ respectively (see again Figure \ref{fig:fusion}). More precisely, there exist $h_q \in [i_q+1..j_q]$ and $c_q \in [a_q+1..b_q]$ such that $\pi_1$ (resp. $\pi_2$) has an occurrence in the intervals $[i_q..h_q-1]$ (resp. $[h_q..j_q]$) of indices and $[a_q..c_q-1]$ (resp. $[c_q..b_q]$) of values. This argument holding for all $q \in [1..K]$, it becomes clear that $\pi_1$ (resp. $\pi_2$) is a common separable pattern between $\sigma^1, \ldots, \sigma^K$ that uses the intervals $[i_q..h_q-1]$ (resp. $[h_q..j_q]$) of indices and $[a_q..c_q-1]$ (resp. $[c_q..b_q]$) of values in $\sigma^q$, for all $q \in [1..K]$. What remains to prove is that $\pi_1$ and $\pi_2$ are of maximal length among all such patterns.

Assume that $\pi_1$ is \emph{not} of maximal length among the common separable patterns between $\sigma^1, \ldots, \sigma^K$ using interval $[i_q..h_q-1]$ of indices and interval $[a_q..c_q-1]$ of values in $\sigma^q$, for all $q \in [1..K]$. Then, there exists $\pi_1^{long}$, a common separable pattern between $\sigma^1, \ldots, \sigma^K$ using interval $[i_q..h_q-1]$ of indices and interval $[a_q..c_q-1]$ of values in $\sigma^q$, for all $q \in [1..K]$, such that $|\pi_1^{long}| > |\pi_1|$. Now by Lemma \ref{lem:fusion}, $\pi_1^{long} \oplus \pi_2$ is a common separable pattern between $\sigma^1, \ldots, \sigma^K$ using the intervals $[i_q..j_q]$ of indices and $[a_q..b_q]$ of values in $\sigma^q$, for all $q \in [1..K]$. And we have $|\pi_1^{long} \oplus \pi_2| > |\pi_1 \oplus \pi_2| = |\pi|$, contradicting the maximality of $\pi$. So $\pi_1$ is a common separable pattern of maximal length between $\sigma^1, \ldots, \sigma^K$ among those using interval $[i_q..h_q-1]$ of indices and interval $[a_q..c_q-1]$ of values in $\sigma^q$, for all $q \in [1..K]$. In the same way, we prove that $\pi_2$ is a common separable pattern of maximal length between $\sigma^1, \ldots, \sigma^K$ among those using interval $[h_q..j_q]$ of indices and interval $[c_q..b_q]$ of values in $\sigma^q$, for all $q \in [1..K]$, ending the proof of the lemma.
\qed
\end{proof}

\begin{proposition}
Algorithm \ref{alg:separable} is correct: it outputs a longest common separable pattern between the $K$ permutations in input. \label{prop:correctness}
\end{proposition}

\begin{proof}
Consider the array $M$ returned by Algorithm \ref{alg:separable}. We show by induction on $\sum_q (j_q-i_q) +(b_q-a_q)$ that $M(i_1,j_1,a_1,b_1, \ldots ,i_K,j_K,a_K,b_K)$ contains a common separable pattern $\pi$ between $\sigma^1, \ldots, \sigma^K$ that is of maximal length among those using, for any $q \in [1..K]$, intervals $[i_q..j_q]$ of indices and $[a_q..b_q]$ of values in $\sigma^q$.

First, there is no loss of generality in assuming that $i_q \leq j_q$ and $a_q\leq b_q$ for all $q \in [1..K]$, since otherwise the above statement is clearly true (indeed, $M(i_1,j_1,a_1,b_1, \ldots ,i_K,j_K,a_K,b_K)$ contains $\epsilon$ and either $[i_q..j_q] = \emptyset$ or $[a_q..b_q] = \emptyset$ for some $q \in [1..K]$).

If $\sum_q (j_q-i_q) +(b_q-a_q) =0$, then $i_q=j_q$ and $a_q=b_q$ for all $q \in [1..K]$. Consequently, the pattern we would like $M(i_1,j_1,a_1,b_1, \ldots ,i_K,j_K,a_K,b_K)$ to contain is a longest common separable pattern between $\sigma^1, \ldots, \sigma^K$ that uses only index $i_q=j_q$ and value $a_q=b_q$ in $\sigma^q$ for all $q \in [1..K]$. Such a pattern is either $\epsilon$ or $1$. And it is $1$ if and only if $\forall q \in [1..K], \sigma^q_{i_q} = a_q$ that is to say if and only if $\forall q \in [1..K], \exists h_q \in [i_q..j_q]$ such that $\sigma^q_{h_q} \in [a_q..b_q]$. On lines 8 to 12 of Algorithm \ref{alg:separable}, we see that in this case $M(i_1,j_1,a_1,b_1, \ldots ,i_K,j_K,a_K,b_K)$ is set correctly.

If $\sum_q (j_q-i_q) +(b_q-a_q) > 0$, we must consider two subcases:

If $\exists q \in [1..K]$ such that $i_q=j_q$ or $a_q=b_q$, let us call $\pi$ a common separable pattern between $\sigma^1, \ldots, \sigma^K$ that is of maximal length among those using, for any $q \in [1..K]$, intervals $[i_q..j_q]$ of indices and $[a_q..b_q]$ of values in $\sigma^q$. Then, just as before, $\pi$ is either $1$ or $\epsilon$, because in at least one of the permutations, say $\sigma^q$, the occurrence of $\pi$ can use at most one index (if $i_q=j_q$) or at most one value (if $a_q=b_q$). And again, lines 8 to 12 of Algorithm \ref{alg:separable} show that $M(i_1,j_1,a_1,b_1, \ldots ,i_K,j_K,a_K,b_K)$ contains $1$ or $\epsilon$. More precisely, the condition on line 9 ensures that it contains $1$ exactly when $1$ has an occurrence in $\sigma^q$ in the intervals $[i_q..j_q]$ of indices and $[a_q..b_q]$ of values, for any $q \in [1..K]$.

It remains to consider the recursive case when $\forall q \in [1..K], i_q < j_q$ and $a_q<b_q$. In this case, consider $\pi$ a common separable pattern $\pi$ between $\sigma^1, \ldots, \sigma^K$ that is of maximal length among those such that, for any $q \in [1..K]$, $\pi$ has an occurrence in $\sigma^q$ in the intervals $[i_q..j_q]$ of indices and $[a_q..b_q]$ of values. Since $\pi$ is separable, then either $\pi = \epsilon$, or $\pi=1$, or $\pi= \pi_1 \oplus \pi_2$, or $\pi= \pi_1 \ominus \pi_2$ where $\pi_1$ and $\pi_2$ are smaller but non-empty separable patterns. Now consider $\pi_{algo}$ a longest pattern in the set $S_{\oplus} \cup S_{\ominus} \cup S $ where 
\begin{eqnarray*}
 S_{\oplus} & = & \{ M(i_1,h_1-1,a_1,c_1-1, \ldots ,i_K,h_K -1,a_K,c_K-1)  \oplus  M(h_1,j_1,c_1,b_1, \\ & & \ldots ,h_K,j_K,c_K,b_K) : i_q < h_q \leq j_q , a_q < c_q \leq b_q , \forall q \in [1..K] \} \\
 S_{\ominus} & = & \{ M(i_1,h_1-1,c_1,b_1, \ldots ,i_K,h_K -1,c_K,b_K)  \ominus  M(h_1,j_1,a_1,c_1-1, \\ & & \ldots ,h_K,j_K,a_K,c_K-1) : i_q < h_q \leq j_q , a_q < c_q \leq b_q , \forall q \in [1..K] \} \\
S_{~} & = & \{1\} \textrm{ if } \forall q \in [1..K], \exists h_q \in [i_q..j_q] \textrm{ such that } \sigma^q_{h_q} \in [a_q..b_q],\\
 & = &  \{\epsilon\} \textrm{ otherwise.}
\end{eqnarray*}
By induction hypothesis, each entry of $M$ that appears in the set $S_{\oplus} \cup S_{\ominus}$ is a common separable pattern between $\sigma^1, \ldots, \sigma^K$ whose occurrence in $\sigma^q$, for any $q \in [1..K]$, uses indices and values in the prescribed intervals, and that is of maximal length among all such patterns. 
Notice also that $S=1$ if and only if $1$ has an occurrence in $\sigma^q$ in the intervals $[i_q..j_q]$ of indices and $[a_q..b_q]$ of values, for all $q \in [1..K]$.
An immediate consequence of those two facts and of Lemma \ref{lem:fusion} is that $\pi_{algo}$ is a common separable pattern between $\sigma^1, \ldots, \sigma^K$ which has, for any $q \in [1..K]$, an occurrence in $\sigma^q$ in the intervals $[i_q..j_q]$ of indices and $[a_q..b_q]$ of values. What is left is to prove  that $|\pi_{algo}| = |\pi|$. This is clear when $\pi=\epsilon$ or $1$. So assume that $\pi=\pi_1 \oplus \pi_2$, the case $\pi=\pi_1 \ominus \pi_2$ being very similar. Since $\pi=\pi_1 \oplus \pi_2$ has an occurrence in each $\sigma^q$ in the intervals $[i_q..j_q]$ of indices and $[a_q..b_q]$ of values, by Lemma \ref{lem:cutting}, there exists indices $(h_q)_{q \in [1..K]}$ and values $(c_q)_{q \in [1..K]}$, with $i_q < h_q \leq j_q , a_q < c_q \leq b_q , \forall q \in [1..K]$, such that $\pi_1$ has an occurrence in each $\sigma^q$ in the intervals $[i_q..h_q-1]$ of indices and $[a_q..c_q-1]$ of values and $\pi_2$ has an occurrence in each $\sigma^q$ in the intervals $[h_q..j_q]$ of indices and $[c_q..b_q]$ of values. Lemma \ref{lem:cutting} also states that $\pi_1$ and $\pi_2$ are of maximal length among the common separable patterns in the given intervals of indices and values. So by induction hypothesis, $|M(i_1, h_1-1, a_1, c_1-1, \ldots , i_K, h_K-1, a_K, c_K-1)| = |\pi_1|$ and $|M(h_1, j_1, c_1, b_1, \ldots , h_K, j_K, c_K, b_K)| = |\pi_2|$. Consequently, $|\pi| = |\pi_1 \oplus \pi_2| = |M(i_1, h_1-1, a_1, c_1-1, \ldots , i_K, h_K-1, a_K, c_K-1)| + |M(h_1, j_1, c_1, b_1, \ldots , h_K, j_K, c_K, b_K)| \leq |\pi_{algo}|$. The inequality $|\pi| \geq |\pi_{algo}|$ being obvious by maximality of $\pi$, we conclude that $|\pi| = |\pi_{algo}|$. This ends the proof in the case $\pi = \pi_1 \oplus \pi_2$. For the case $\pi = \pi_1 \ominus \pi_2$, the proof follows the exact same steps, with $\pi_1$ having an occurrence in each $\sigma^q$ in the intervals $[i_q..h_q-1]$ of indices and $[c_q..b_q]$ of values and $\pi_2$ having an occurrence in each $\sigma^q$ in the intervals $[h_q..j_q]$ of indices and 
$[a_q..c_q-1]$ of values.
\qed
\end{proof}

\begin{proposition}
Algorithm \ref{alg:separable} runs in time $\mathcal{O}(n^{6K+1})$ and space $\mathcal{O}(n^{4K+1})$. \label{prop:complexity}
\end{proposition}

\begin{proof}
Algorithm \ref{alg:separable} handles an array $M$ of size $\mathcal{O}(n^{4K})$, where each cell contains a pattern of length at most $n$, so that the total space complexity is $\mathcal{O}(n^{4K+1})$. For filling in one entry $M(i_1,j_1,a_1,b_1, \ldots , i_K, j_K, a_K, b_K)$, if $\exists q \in [1..K]$ such that $i_q=j_q$ or $a_q=b_q$ (lines 9 to 13 of Algorithm 1), the time complexity is $\mathcal{O}(n^K)$. If no such $q$ exists (line 15 of Algorithm 1), the time complexity needed to fill in $M(i_1,j_1,a_1,b_1, \ldots , i_K, j_K, a_K, b_K)$, using the entries of $M$ previously computed, is $\mathcal{O}(n^{2K+1})$. Indeed, we search for an element of maximal length among $\mathcal{O}(n^{2K})$ elements, each element being computed in $\mathcal{O}(n)$-time as the concatenation of two previously computed entries of $M$. Consequently, the total time complexity to fill in $M$ is $\mathcal{O}(n^{6K+1})$.
\qed
\end{proof}

A consequence of propositions \ref{prop:correctness} and \ref{prop:complexity} is:
\begin{theorem}
For any fixed integer $K$, the problem of computing a longest common separable pattern between $K$ permutations is in $P$.
\end{theorem}

%Notice that when the number of permutations in input is finite but not known in advance, our algorithm still apply, but is not polynomial. Indeed, we show in the next section that the problem of computing a longest common separable pattern between a finite but arbitrary number of permutations is $NP$-hard.\\

We may wonder whether a longest common separable pattern between two permutations $\sigma^1$ and $\sigma^2$ (computed in polynomial time by Algorithm \ref{alg:separable}) is a good approximation of a longest common pattern between $\sigma^1$ and $\sigma^2$ (whose computation is $NP$-hard). Section \ref{section:approximation} gives a negative answer to this question, by the more general Corollary \ref{cor:sqrt}.

\section{Hardness result}
\label{section:hardness result}

We proved in the preceding section that the \PB{\LCSP} problem is
polynomial-time solvable provided a constant number of input
permutations. 
We show here that the \PB{\CSP} problem (the general decision version of \PB{\LCSP}), is \NPC.

\begin{problem}
[\PB{\CSP} Problem]~\\
\textsc{Input}: A set $X$ of permutations and an integer $k$.\\
\textsc{Output}: A boolean indicating if there a separable pattern of length $k$ occurring in each $\sigma \in X$.
\label{problem:CSP}
\end{problem}

Actually, we will prove more, namely that the
\PB{\CSP} problem is \NPC even if each input permutation is separable. 
An immediate consequence is the $NP$-hardness of \PB{\LCSP}.
For ease of exposition, our proof is given in terms of matching
diagrams. 

\begin{definition}[Matching Diagram]
A \emph{matching diagram} $G$ of size $n$ 
is a vertex-labeled graph of order \textit{i.e.,} number of vertices, $2n$ and size \textit{i.e.,} number of edges, $n$ where each vertex
is labeled by a distinct label from $\{1, 2, \ldots, 2n\}$ and
each vertex $i \in \{1, 2, \ldots, n\}$ 
(resp. $j \in \{n+1, n+2, \ldots, 2n\}$)
is connected by an edge to exactly one vertex 
$j \in \{n+1, n+2, \ldots, 2n\}$ 
(resp. $i \in \{1, 2, \ldots, n\}$). 
We denote the set of vertices and edges of $G$ by $V(G)$ and $E(G)$,
respectively. 
\end{definition}

\begin{figure}
\begin{center}
%% \begin{tikzpicture}[scale=.4]
%% \draw (-2,3) node {$G$:\;\;};
%% %
%% \draw[rounded corners] 
%% (0,3.2) .. controls (0,7) and (27,7) .. (27,3.2);
%% \draw[rounded corners] 
%% (3,3.2) .. controls (3,6) and (15,6) .. (15,3.2);
%% \draw[rounded corners] 
%% (6,3.2) .. controls (6,6) and (18,6) .. (18,3.2);
%% \draw[rounded corners] 
%% (9,3.2) .. controls (9,6) and (24,6) .. (24,3.2);
%% \draw[rounded corners] 
%% (12,3.2) .. controls (12,5) and (21,5) .. (21,3.2);
%% %
%% \foreach \x in {0,3,6,9,12,15,18,21,24,27} 
%% \shadedraw (\x,3) circle (2ex);
%% %
%% \draw (0,2) node {$1$};
%% \draw (3,2) node {$2$};
%% \draw (6,2) node {$3$};
%% \draw (9,2) node {$4$};
%% \draw (12,2) node {$5$};
%% \draw (15,2) node {$6$};
%% \draw (18,2) node {$7$};
%% \draw (21,2) node {$8$};
%% \draw (24,2) node {$9$};
%% \draw (27,2) node {$10$};
%% %
%% \draw (13,0) node[red] {$\pi = $};
%% \draw (15,0) node[red] {$2$};
%% \draw (18,0) node[red] {$3$};
%% \draw (21,0) node[red] {$5$};
%% \draw (24,0) node[red] {$4$};
%% \draw (27,0) node[red] {$1$};
%% \end{tikzpicture}
\includegraphics{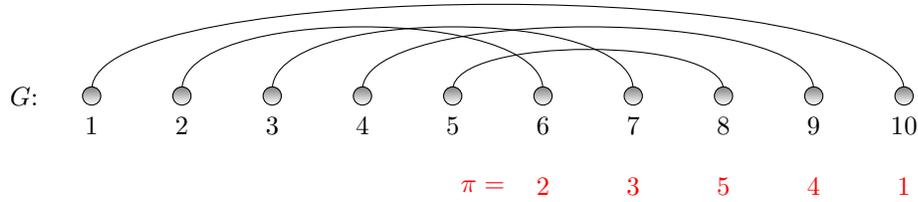}
\caption{\label{fig:bijection}Shown here is the correspondence between the
  permutation $\pi = 2\;3\;5\;4\;1$ and the associated matching
  diagram $G$.}
\end{center}
\end{figure}

It is well-known that matching diagrams of size $n$ are in one-to-one 
correspondence with permutations of length $n$ 
(see Figure~\ref{fig:bijection} for an illustration).

Let $G$ and $G'$ be two matching diagrams.
The matching diagram $G'$ is said to \emph{occur} in $G$ 
if one can obtain $G'$ from $G$ by a sequence of edge deletions. 
More formally, the deletion of the edge $(i,j)$, $i < j$, consists in 
(1) the deletion of the edge $(i,j)$,
(2) the deletion of the two vertices $i$ and $j$, and 
(3) the relabeling of all vertices $k \in [i+1..j-1]$ to $k-1$ and all the vertices $k>j$ to $k-2$. Therefore, the decision version of the \PB{\LCSP} is equivalent to the following problem: Given a set of matching diagrams and a positive integer $k$, find a matching diagram of size $k$ which occurs in each input diagram \cite{CPM2006}.

Clearly, two edges in a matching diagram $G$ are either crossing 
%% \begin{tikzpicture}[scale=.2]
%% \draw[rounded corners]
%% (0,0) .. controls (1,1) and (2,1) .. (3,0);
%% \draw[rounded corners]
%% (1,0) .. controls (2,1) and (3,1) .. (4,0);
%% \foreach \x in {0,1,3,4} 
%% \shadedraw (\x,0) circle (1ex);
%% \end{tikzpicture}
\includegraphics{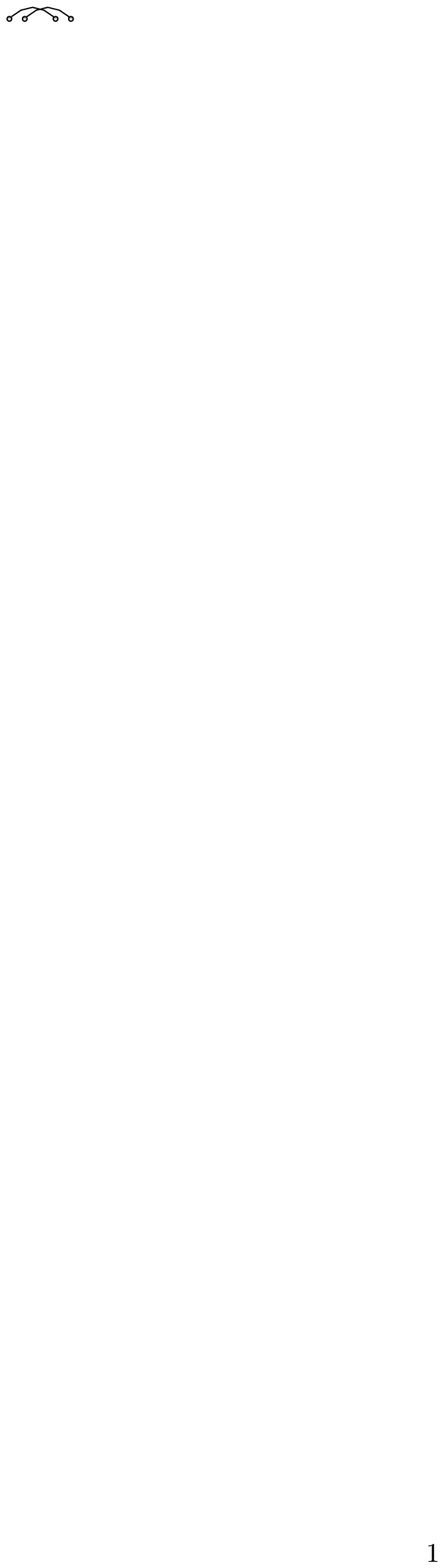}
or
nested
\includegraphics{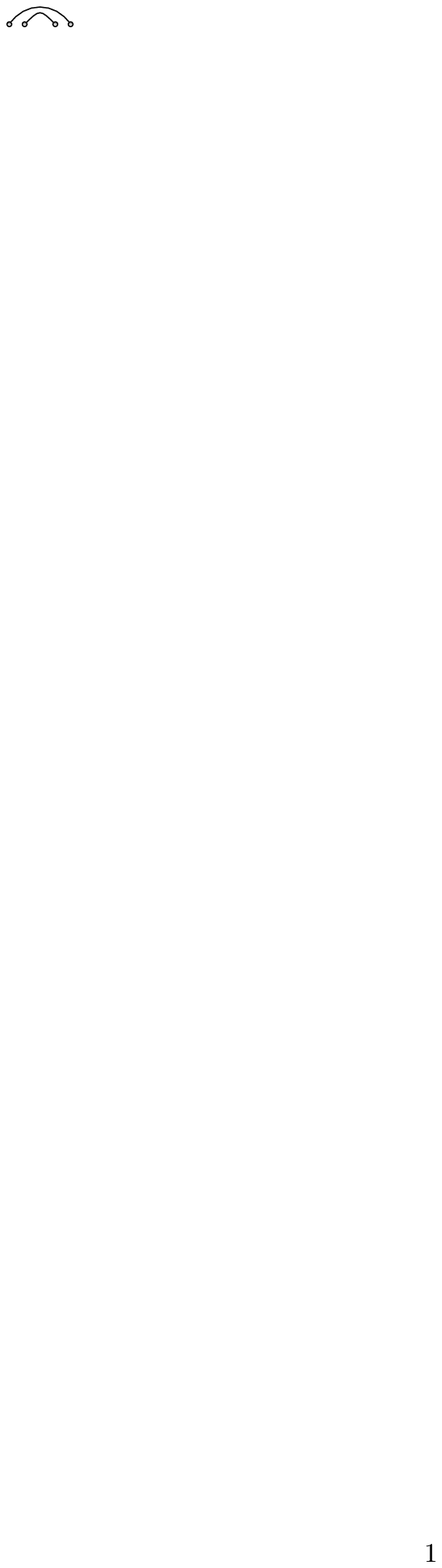}
%% \begin{tikzpicture}[scale=.2]
%% \draw[rounded corners]
%% (1,0) .. controls (2,1) .. (3,0);
%% \draw[rounded corners]
%% (0,0) .. controls (1,1.5) and (3,1.5) .. (4,0);
%% \foreach \x in {0,1,3,4} 
%% \shadedraw (\x,0) circle (1ex);
%% \end{tikzpicture}
.
Moreover, it is easily seen that an occurrence in $G$ of a matching
diagram $G'$ of which all edges are crossing (resp. nested) correspond
to an occurrence in the permutation associated with $G$ of
an \emph{increasing} (resp. \emph{decreasing}) subsequence.

For the purpose of permutations, convenient matching diagrams are
needed.
A matching diagram is called a \emph{tower} if it is composed of
pairwise nested edges
%% \begin{tikzpicture}[scale=.15]
%% \draw[rounded corners]
%% (1,0) .. controls (2,1) .. (3,0);
%% \draw[rounded corners]
%% (0,0) .. controls (1,1.5) and (3,1.5) .. (4,0);
%% \draw[rounded corners]
%% (-1,0) .. controls (0,2) and (4,2) .. (5,0);
%% \foreach \x in {-1,0,1,3,4,5} 
%% \shadedraw (\x,0) circle (1ex);
%% \end{tikzpicture}
\includegraphics{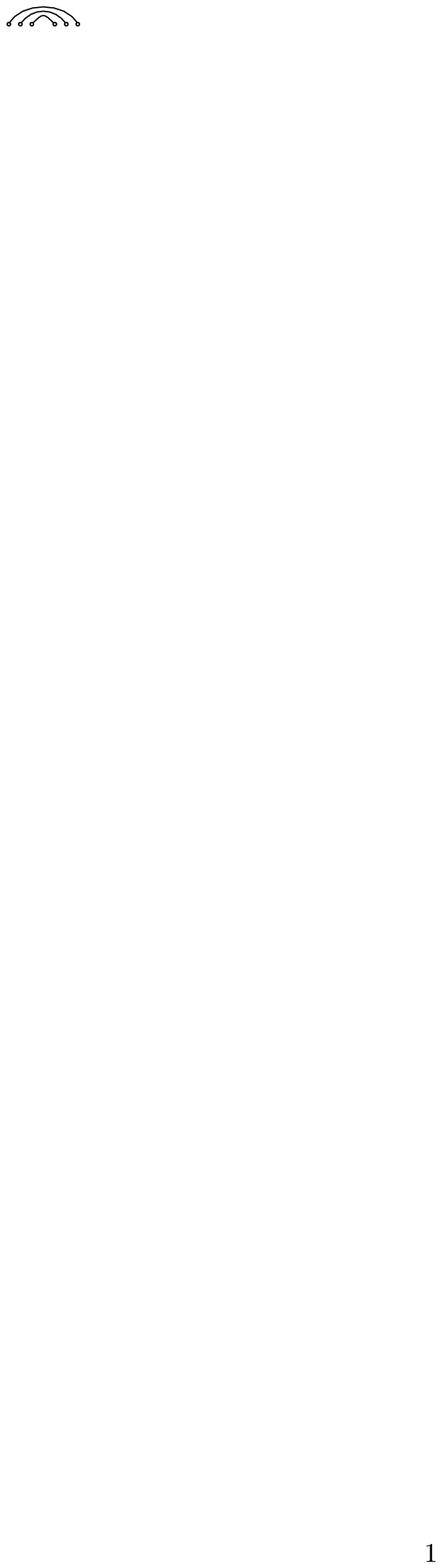}
and a \emph{staircase} if it is
composed of pairwise crossing edges
%% \begin{tikzpicture}[xscale=.15,yscale=0.2]
%% \draw[rounded corners]
%% (0,0) .. controls (1,1) and (2,1) .. (3,0);
%% \draw[rounded corners]
%% (1,0) .. controls (2,1) and (3,1) .. (4,0);
%% \draw[rounded corners]
%% (2,0) .. controls (3,1) and (4,1) .. (5,0);
%% \foreach \x in {0,1,2,3,4,5} 
%% \shadedraw (\x,0) circle (1ex);
%% \end{tikzpicture}
\includegraphics{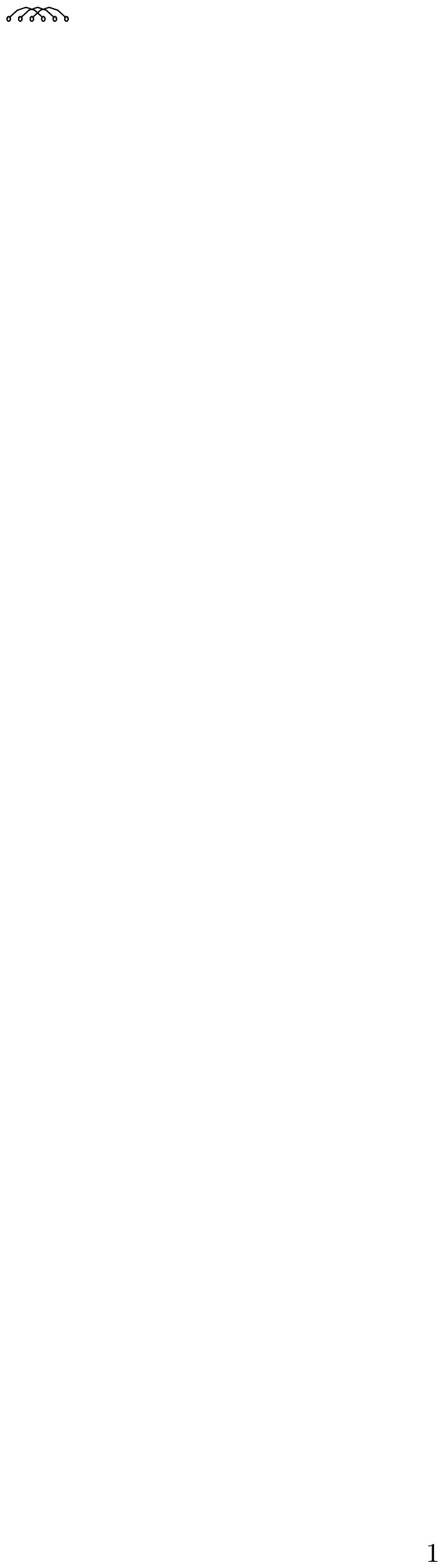}
.
A matching diagram is called a \emph{tower of staircases} if its edge
set can be partitioned in nested staircases
%% \begin{tikzpicture}[xscale=.15,yscale=0.1]
%% \draw[rounded corners]
%% (0,0) .. controls (2,4) and (6,4) .. (8,0);
%% \draw[rounded corners]
%% (1,0) .. controls (3,4) and (7,4) .. (9,0);
%% \draw[rounded corners]
%% (3,0) .. controls (4,1.5) .. (5,0);
%% \draw[rounded corners]
%% (4,0) .. controls (5,1.5) .. (6,0);
%% \foreach \x in {0,1,3,4,5,6,8,9} 
%% \shadedraw (\x,0) circle (1ex);
%% \end{tikzpicture}
\includegraphics{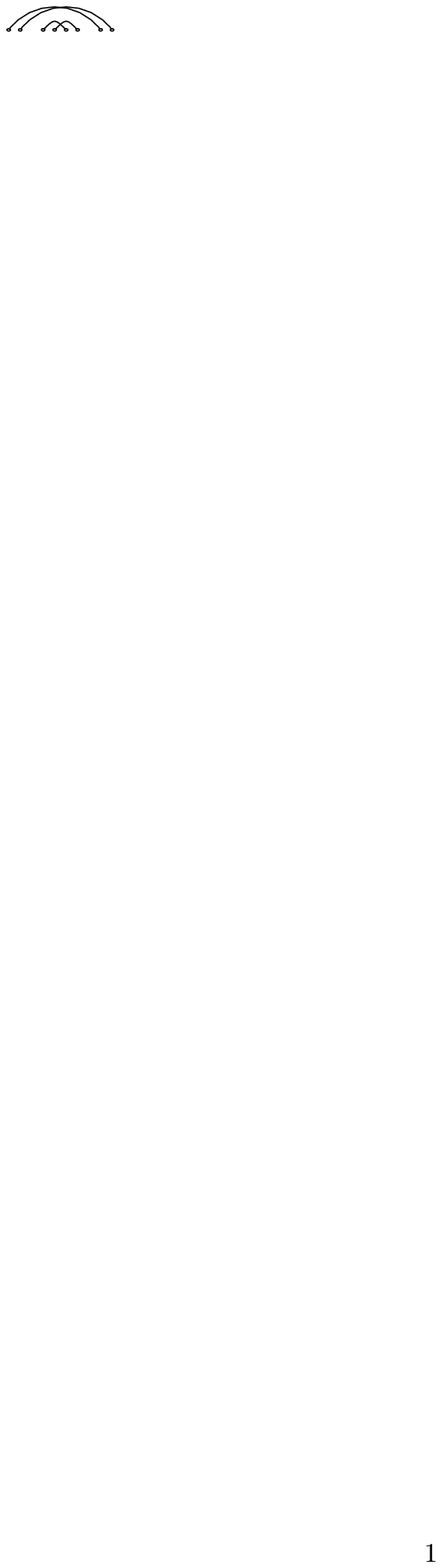}.

\begin{figure}
\begin{center}
\includegraphics{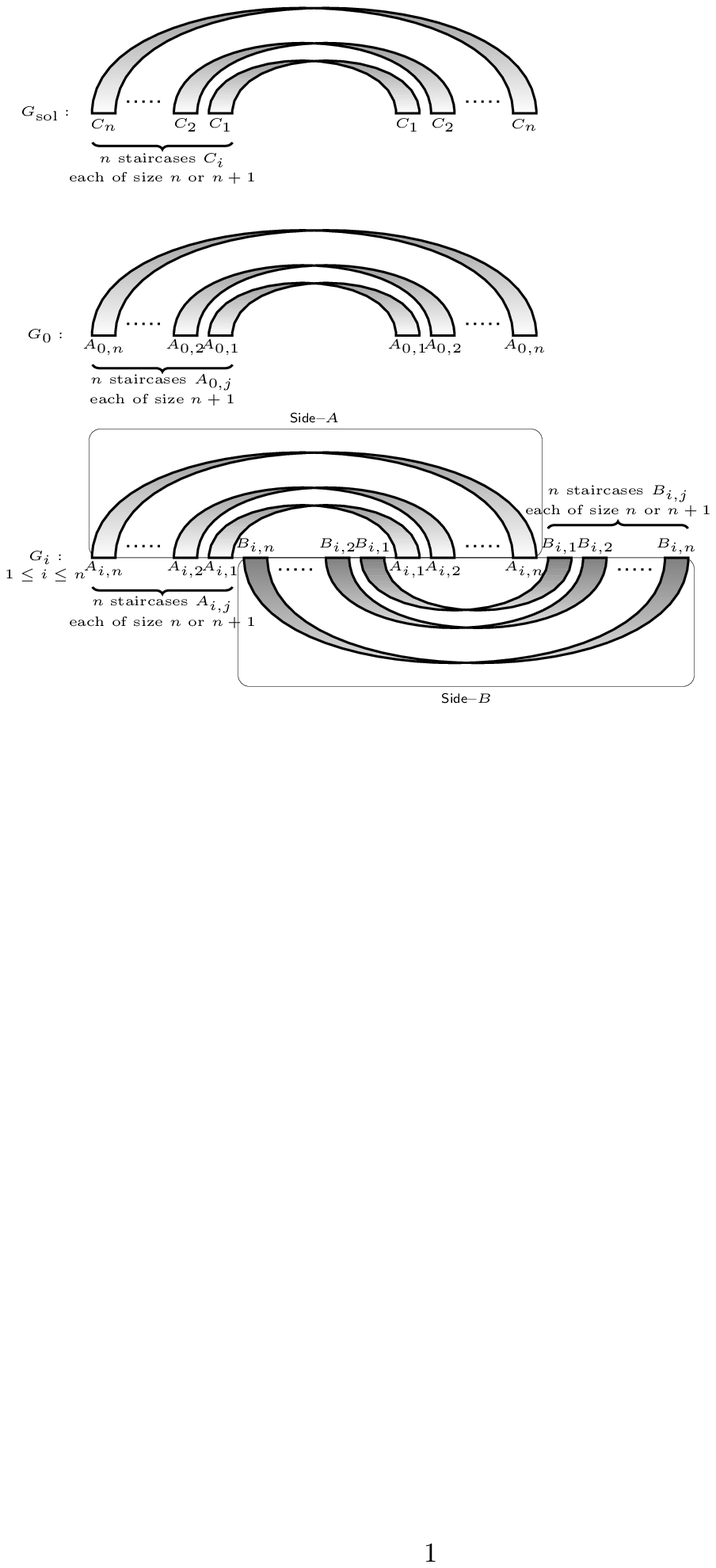}
\caption{\label{fig:NP-hard}Reduction in the proof of Theorem \ref{thm:hardness}}
\end{center}
\end{figure}

\begin{theorem}
The \PB{\CSP} problem is \NPC even if each input permutation is
separable. \label{thm:hardness}
\end{theorem}
\begin{proof}
\PB{\CSP} is clearly in $NP$. For proving hardness, we reduce from the \PB{Independent-Set} problem which is known to be 
\NPC \cite{garey79} .
Let $G$ be an arbitrary graph instance of the \PB{Independent-Set}
problem.
Write $V(G) = \{1, 2, \ldots, n\}$.
We now detail the construction of $n+1$ matching diagrams
$G_0, G_1,$ $ G_2, \ldots, G_n$, each corresponding to a separable
permutation.
First the matching diagram $G_0$ is a tower of $n$ staircases 
$A_{0,1}, A_{0,2}, \ldots, A_{0,n}$, each of size $n+1$
(see Figure~\ref{fig:NP-hard}, middle part; staircases are represented
by shaded forms), \emph{i.e.},
$$
\forall j,\; 1 \leq j \leq n,
\qquad
|A_{0,j}| = n+1\text{.}
$$
Each matching diagram $G_i$, $1 \leq i \leq n$, is composed of two
crossing towers of $n$ staircases each referred as 
$A_{i,1}, A_{i,2}, \ldots, A_{i,n}$ and
$B_{i,1}, B_{i,2}, \ldots, B_{i,n}$
(see Figure~\ref{fig:NP-hard}, bottom part),
and defined as follows:
\begin{alignat*}{3}
&\forall i,\; 1 \leq i \leq n,
\;\;
\forall j,\; 1 \leq j \leq n,
&\qquad
|A_{i,j}| &\begin{cases}
n+1 & \text{if $i \neq j$} \\
n & \text{if $i = j$}
\end{cases}
\\
&\forall i,\; 1 \leq i \leq n,
\;\;
\forall j,\; 1 \leq j \leq n,
&\qquad
|B_{i,j}| &\begin{cases}
n+1 & \text{if $(i,j) \notin E(G)$} \\
n   & \text{if $(i,j) \in E(G)$.}
\end{cases}
\end{alignat*}
It is simple matter to check that all matching diagrams $G_i$,
$0 \leq i \leq n$, correspond to separable permutations
and that our construction can be carried on in polynomial time. 
This ends our construction. 

We claim that there exists an independent set $V' \subseteq V(G)$ of
size $k$ in $G$
if and only if
there exists a matching diagram $\Gsol$ of
size $n^2+k$ that occurs in each input matching diagram $G_i$, 
$0 \leq i \leq n$. 

Suppose that there exists an independent set $V' \subseteq V(G)$ of
size $k$ in $G$.
Consider as a solution $\Gsol$ the tower of $n$ staircases 
$C_1, C_2, \ldots, C_n$ (see Figure~\ref{fig:NP-hard}, top part) of
total size $n^2 + k$, where the size of each staircase is defined
according to the following: 
$|C_i| = n$ if $i \notin V'$ and 
$|C_i| = n+1$ if $i \in V'$. 
We claim that $\Gsol$ occurs in each input matching diagram
$G_i$, $0 \leq i \leq n$.
Indeed, by construction,
for any $i \in V(G)$,
if $i \notin V'$ then $\Gsol$ occurs in 
\textsf{Side}--$A$ of $G_i$ and
if $i \in V'$ then $\Gsol$ occurs in 
\textsf{Side}--$B$ of $G_i$.
To complete the claim,
we note that $\Gsol$ occurs in $G_0$ 
($G_0$ is indeed a tower of $n$ staircases, each of size $n+1$). 

Conversely, suppose that there exists a matching diagram $\Gsol$ of
size $n^2+k$ that occurs in each input matching diagram $G_i$, 
$0 \leq i \leq n$. 
Let us prove the following:
\begin{claim}
$\Gsol$ is a tower of $n$ staircases. Furthermore, for any $i \in [1..n]$, $\Gsol$ occurs completely in \textsf{Side}--$A$ or completely in \textsf{Side}--$B$ in $G_i$. \label{claim}
\end{claim}
\begin{proof}
Let us first focus on an occurrence of $\Gsol$ in $G_0$.
Since $G_0$ is a tower of staircases, then it follows that
$\Gsol$ is a tower of staircases as well.
We now prove that $\Gsol$ is a tower of $n$ staircases.
Indeed, suppose, aiming at a contradiction, that 
$\Gsol$ is a tower of $n' < n$ staircases.
Then it follows that $\Gsol$ has size at most
$(n+1)(n-1) = n^2 - 1 < n^2 + k$.
This is the desired contradiction, and hence $\Gsol$ is a tower of $n$
staircases $C_1, C_2, \ldots, C_n$, each of size at most $n+1$. 
We now turn to considering an occurrence of $\Gsol$ in $G_i$,
$1 \leq i \leq n$.
We prove that $\Gsol$ occurs completely in \textsf{Side}--$A$ or 
completely in \textsf{Side}--$B$ (see Figure~\ref{fig:NP-hard}).
Suppose, for the sake of contradiction, 
that $\Gsol$ matches at least one edge in \textsf{Side}--$A$, say $e$,
and at least one edge in \textsf{Side}--$B$, say $e'$.
By construction, $e$ and $e'$ are crossing edges in $G_i$, and hence 
$e$ and $e'$ are matched by two edges that belong to the same
staircase in $\Gsol$ 
($\Gsol$ is indeed a tower of staircases). 
We now observe that
any edge in \textsf{Side}--$A$ crosses the edge $e'$ and
any edge in \textsf{Side}--$B$ crosses the edge $e$.
Then it follows that the occurrence of $\Gsol$ in $G_i$ induces a single 
staircase. 
But $\Gsol$ is a tower of $n$ staircases, each of size at most $n+1$. 
A contradiction.
Therefore, in any matching diagram $G_i$, $1 \leq i \leq n$,
$\Gsol$ occurs completely in \textsf{Side}--$A$ or 
completely in \textsf{Side}--$B$.
\end{proof}
As an important consequence of the claim, 
there is thus no loss of generality in assuming that each staircase of  
$\Gsol$ has size $n$ or $n+1$, and hence $\Gsol$ is composed of $n-k$
staircases of size $n$ and $k$ staircases of size $n+1$ 
(since $\Gsol$ has size $n^2 + k$).
Consider now the subset $V' \subseteq V(G)$ defined as follows: 
$i \in V(G)$ belongs to $V'$ if and only if
$|C_i| = n+1$, \emph{i.e.}, the staircase $C_i$ has size $n+1$ in
$\Gsol$.
According to the above, $V'$ has certainly size $k$.  
We shall show that $V'$ is an independent set in
$G$. 
Indeed, let $i \in V'$ and consider the matching diagram $G_i$.
Since $i \in V'$ then it follows that $|C_i| = n+1$, and hence
$\Gsol$ occurs in \textsf{Side}--$B$ in $G_i$
(the latter follows from the fact that $|A_{i,i}| = n$).
But, by construction, for $1 \leq j \leq n$, $|B_{i,j}| = n+1$ if and only
if $(i,j)$ is not an edge in $G$ (in particular, $|B_{i,i}| = n+1$).
Hence, the vertex $i$ is not connected to any vertex in $V'$.
Therefore, since the argument applies to each matching diagram $G_i$
with $i \in V'$, $V'$ is an independent set in $G$.
\qed
\end{proof}

\section{Approximation ratio} 
\label{section:approximation}

In this section, we return to the \PB{\LCCP} Problem for $K$ permutations. As said before, the general \PB{\LCP} Problem is NP-hard as well as the Pattern Involvement Problem \cite{BBL98}. In this section we prove the following result:

\begin{theorem}\label{thm:racine}
For all $\epsilon > 0$ and \C, a pattern-avoiding class of permutations, there exists a sequence $(\sigma_n)_{n \in {\mathbb N}}$ of permutations $\sigma_n \in S_n$ such that 
$$|\pi_n| = o\left(n^{0.5+\epsilon}\right)$$

where $\pi_n$ is the longest pattern of class \C involved in $\sigma_n$.
\end{theorem}

Before proving this result we need the following Lemma.
\begin{lemma}\label{lem:comptagePattern}
Given a permutation $\pi \in S_k$, the number of permutations $\sigma \in S_n$ such that $\pi$ is involved in $\sigma$ is at most $(n-k)!\binom{n}{k}^2$.
\end{lemma}
\begin{proof}
Let $\pi=\pi_1\pi_2\ldots,\pi_k$ and $\sigma=\sigma_1\sigma_2\ldots\sigma_n$ be such that $\pi$ is involved in $\sigma$. There exist $i_1 < i_2< \ldots < i_k$ such that $\sigma_{i_1}\sigma_{i_2}\ldots\sigma_{i_k}$ is order-isomorphic to $\pi$. 
Then $\sigma = u_1\sigma_{i_1}u_2\sigma_{i_2}\ldots u_k\sigma_{i_k}u_{k+1}$ where $u_i$ is a factor of $\sigma$. We will call $u_i$'s the blocks associated to $\sigma, \pi$ and an occurrence of $\pi$ in $\sigma$.

Conversely, suppose we are given a permutation $\pi \in S_k$, and $u_1,u_2,\ldots ,u_{k+1}$ $k+1$ sequences of distinct numbers in $\{1\ldots n\}$ such that $|u_1|+|u_2|+\ldots +|u_{k+1}| = n-k$. Let $E=\{n_1,n_2,\ldots n_k\}$ the $k$ integers in $\{1\ldots n\}$ which do not appear in a block $u_i$. We denote by $\bar j$ the $j^{th}$ smallest element in $E$.
Then the permutation $\sigma= u_1\overline{\pi_1}u_2\overline{\pi_2}\ldots\overline{\pi_k}u_{k+1}$  is a permutation of $S_n$ and $\pi$ is involved in $\sigma$. 
For example, if $\pi = 2143$, $n=9$, $u_1=31$, $u_2=\varnothing$, $u_3=8$, $u_4=65$ and $u_5=\varnothing$ then $E=\{2,4,7,9\}$ and $\sigma=31\textbf{4}\textbf{2}8\textbf{9}65\textbf{7}$.

Note that two different lists of blocks could lead to the same permutation $\sigma$ if $\pi$ has several occurrences in $\sigma$ as shown in Figure \ref{fig:severalPatternsOccurence}. 

\begin{figure}[ht]
\begin{center}
$$
2431 \longrightarrow
\begin{cases}
 (1,\varnothing, \varnothing, 56, \varnothing)\\
 (1,\varnothing, 4, 6, \varnothing)\\
 (1,\varnothing, 45, \varnothing, \varnothing)\\
\end{cases}
\longrightarrow 1374562$$
\caption{Relation between pattern, blocks and permutation}
\label{fig:severalPatternsOccurence}
\end{center}
\end{figure}
Thus, the number of permutations $\sigma$ such that $\pi$ is involved in $\sigma$ is {\it at most} the number of different lists of blocks. There are $\binom{n}{n-k}$ different choices for the numbers that appear in one of the $u_i$. Then there are $(n-k)!$ different orders for these numbers. The last step is to cut the word so obtained into $k+1$ (possibly empty) blocks $u_1, u_2, \ldots u_{k+1}$. There are $\binom{n}{n-k}$ such choices. Hence we have the claimed formula.
\qed
\end{proof}

We can now prove Theorem \ref{thm:racine}.

\begin{proof}
We make the proof by contradiction. We first prove that if the result were false,  every permutation of length $n$ would contain a pattern of \C of length {\em exactly} $k=\lceil n^{0.5+\epsilon} \rceil$. Next, we show that the number of permutations of length $n$ containing one permutation of $\C \bigcap S_{k}$ as a pattern is strictly less than $n!$.

Suppose that there exist $\epsilon >0$ and \C a pattern-avoiding class of permutations such that for every permutation $\sigma \in S_n$, the longest pattern $\pi \in \C$ of $\sigma$ has length $|\pi| \geq \lceil |\sigma|^{0.5+\epsilon}\rceil =k$. As \C is closed - every pattern $\pi$ of a permutation $\sigma \in \C$ is also in \C - for every permutation $\sigma \in S_n$ there exists a  pattern $\pi \in \C$ of $\sigma$ whose length is {\em exactly} $\ |\pi| = k$. 

But the number of permutations in $\C \bigcap S_k$ is at most $c^{k}$ by \cite{MT04}. By Lemma \ref{lem:comptagePattern}, for each permutation in $\C \bigcap S_k$, the number of permutations in $S_n$ having this permutation as a pattern is at most $(n-k)!\binom{n}{k}^2$ . Thus the number of permutations in $S_n$ having a pattern in $\C \bigcap S_{\geq k}$ is at most $c^k (n-k)!\binom{n}{k}^2$. But with $k=\lceil n^{0.5+\epsilon} \rceil$, $c^k (n-k)!\binom{n}{k}^2 = o\left( n^{n^{1-2\epsilon}} \right) = o\left(n!\right)$. Note that a similar proof is given in \cite{EELW02} for finding the smallest permutation containing all patterns of a given length.
\qed
\end{proof}

\begin{corollary}
The \PB{\LCP} Problem cannot be approximated with a tighter ratio than $\sqrt{Opt}$ by the \PB{\LCCP} Problem, where \C is a pattern-avoiding class of permutations, and $Opt$ is the size of an optimal solution to the \PB{\LCP} Problem. \label{cor:sqrt}
\end{corollary}

\begin{proof}
Consider the \PB{\LCP} Problem between $\sigma$ and $\sigma$. Then the optimal solution to the \PB{\LCP} Problem is $\sigma$. But the solution to the \PB{\LCCP} Problem is a longest pattern of $\sigma$ belonging to the class \C. By Theorem \ref{thm:racine}, such a pattern may have size $\sqrt{|\sigma|}$ asymptotically.
\qed
\end{proof}

%%%%%%%%%%%%%%%%%%%%%%%%%%%%%%%%%%%%%%

\bibliographystyle{plain}
\bibliography{Biblio}

\end{document}